\documentclass[12pt]{article}
\usepackage[mathscr]{eucal}
\usepackage{amssymb,latexsym}
\usepackage{verbatim}
\usepackage{amsmath}
\usepackage{amsthm}
\usepackage{enumerate}
\usepackage{authblk}
\usepackage{color}
\usepackage{url}

\usepackage{tikz}
\usetikzlibrary{matrix}
\usetikzlibrary{snakes}
\usetikzlibrary{arrows,calc,shapes,decorations.pathreplacing}
\usepackage{tikz-cd}

\usepackage[normalem]{ulem}

\usepackage{amscd,amssymb,amsthm,verbatim}
\usepackage{enumerate}
\usepackage{bm}
\usepackage{mathrsfs, graphicx} 
\usepackage{stmaryrd}

\usepackage[nottoc,notlot,notlof]{tocbibind}
\usepackage{setspace}

\usepackage[%
bookmarks=true,
colorlinks,
linkcolor=blue,
urlcolor=blue,
citecolor=blue,
plainpages=false,
pdfpagelabels,
final,
breaklinks=true\between
]{hyperref}
\usepackage{hyperref}

\newtheorem{thm}{Theorem}
\newtheorem{prop}[thm]{Proposition}
\newtheorem{cor}[thm]{Corollary}
\newtheorem{conj}[thm]{Conjecture}

\renewcommand\l{\lambda}

\newcommand\bbR{{\mathbb R}}

\newcommand\wt{\widetilde}

\renewcommand\S{\Sigma}

\newcommand\s{\sigma}
\renewcommand\d{\partial}
\newcommand\e{\varepsilon}
\renewcommand\b{\beta}

\newcommand\<{\langle}
\renewcommand\>{\rangle}

\newcommand\ric{{\rm Ric}}
\renewcommand\l{\lambda}
\newcommand\g{\gamma}

\renewcommand\a{\alpha}

\renewcommand\t{\tau}

\newcommand\bl{\textcolor{blue}}

\newcommand\beq{\begin{equation}}
\newcommand\eeq{\end{equation}}
\newcommand\ben{\begin{enumerate}}
\newcommand\een{\end{enumerate}}
\newcommand\bit{\begin{itemize}}
\newcommand\eit{\end{itemize}}

\DeclareMathOperator{\Ric}{Ric}

\newcommand{\R}{\mathbb R}

\newcommand{\pd}{\partial}

\newcounter{mnotecount}

\title{A CMC existence result for expanding cosmological spacetimes}

\author[1]{Gregory J. Galloway\footnote{galloway@math.miami.edu}}
\author[2]{Eric Ling\footnote{el@math.ku.dk}}

\affil[1]{University of Miami, Coral Gables, FL, USA}

\affil[2]{Copenhagen Centre for Geometry and Topology (GeoTop),
\linebreak
Department of Mathematical Sciences, \linebreak University of Copenhagen, Denmark}

\begin{document}
\date{}
\maketitle

\vspace{.15in}

\begin{abstract} 
We establish a new CMC (constant mean curvature) existence result for cosmological spacetimes, by which we mean globally hyperbolic spacetimes with compact Cauchy surfaces. If a cosmological spacetime satisfying the strong energy condition  contains an expanding Cauchy surface and is future timelike geodesically complete, then the spacetime contains a CMC Cauchy surface.
 This result settles, under certain circumstances, a conjecture of the authors and a conjecture of Dilts and Holst.
Our proof relies on the construction of barriers in the support sense, and the CMC Cauchy surface is found as the asymptotic limit of mean curvature flow. Analogous results are also obtained in the case of a positive cosmological constant $\Lambda > 0$. Lastly, we include some comments concerning the future causal boundary for cosmological spacetimes which pertain to  the CMC conjecture of the authors.

\end{abstract}

\begin{spacing}{.25}
\tableofcontents 
\end{spacing}

\section{Introduction}

As is well known, in order to solve the Einstein equations of general relativity, one must first obtain initial data that satisfies the Einstein constraint equations, and then use the Einstein evolution equations to evolve that data.   
The standard conformal method for solving the constraint equations involves solving certain equations which decouple when the mean curvature associated to the initial data is assumed to be constant.  This assumption then substantially simplifies the problem of solving the constraint equations; see e.g., \cite{Isenberg}.  There are also known advantages for solving the Einstein evolution equations if one works in the CMC gauge, given CMC initial data; see e.g., \cite{Andersson}.  

One is naturally led to consider the circumstances under which a given spacetime  admits a CMC spacelike hypersurface. In \cite{Dilts}, Dilts and Holst address this issue for the class of what they call  {\it cosmological spacetimes}.  These are globally hyperbolic spacetimes $(M,g)$, with compact Cauchy surfaces, that satisfy the strong energy condition, $\ric(X,X) \ge 0$ for all timelike vectors  $X$.   
In this setting, they review various results that establish the existence of CMC Cauchy surfaces under various conditions (e.g.,  as in \cite{Gerhardt, Bart84, Bart88}), as well as examples of spacetimes (suitably global) that fail to have CMC Cauchy surfaces (e.g.,  as in \cite{Bart88, CIP} and generalizations thereof \cite{Ling_Ohanyan1, Ling_Ohanyan2}).

  Consideration of these results and examples led them  to formulate several conjectures concerning the existence (and nonexistence) of CMC Cauchy surfaces in cosmological spacetimes.  Motivated by some of their considerations, in \cite{GalLing} we obtained the following CMC existence result which relies on a certain spacetime curvature condition.

\begin{thm}[\cite{GalLing}]\label{CMCnonpos}
Let $(M,g)$ be a spacetime with compact Cauchy surfaces. Suppose $(M,g)$ is future timelike geodesically complete and has everywhere nonpositive timelike sectional curvatures, i.e., $K \leq 0$ everywhere. Then 
$(M,g)$ contains a CMC Cauchy surface.
\end{thm}

Since the sectional curvature condition above implies the strong energy condition, we were naturally led to the following conjecture.

\begin{conj}\label{conj.}
Let $(M,g)$ be a spacetime with compact Cauchy surfaces. If $(M,g)$ is future timelike geodesically complete and satisfies the strong energy condition, i.e., $\emph{Ric}(X,X) \geq 0$ for all timelike $X$, then $(M,g)$ contains a CMC Cauchy surface.
\end{conj}

In this article, we show that the above conjecture holds if we make the additional assumption that the spacetime contains a Cauchy surface with nonnegative mean curvature. (A similar result was obtained previously under special conditions \cite{GalLing2}.) Specifically, our main result is the following. 

\begin{thm}\label{thm: main}
Let $(M,g)$ be a spacetime with compact Cauchy surfaces, and suppose there exists a smooth spacelike Cauchy surface $V$ with  mean curvature $H \geq 0$. If $(M,g)$ is future timelike geodesically complete and satisfies the strong energy condition, i.e., $\ric(X,X) \geq 0$ for all timelike $X$, then $(M,g)$ contains a CMC Cauchy surface.
\end{thm}

This theorem also settles Conjecture 3.7 in \cite{Dilts} in the case that the spacetime is future timelike geodesically complete. 
(In fact, as discussed in Section \ref{thm: main}, the future completeness assumption can be weakened.)
The proof of 
Theorem \ref{thm: main} 
makes use of certain prescribed mean curvature existence results, via mean curvature flow, of Ecker and Huisken \cite{Ecker_Huisken} modified to apply to the setting of mean curvature barriers in the support sense (as defined in section \ref{thm: main}).  

In Section 2 we provide the relevant background material  concerning mean curvature flow needed for our purposes. The proof of Theorem \ref{thm: main} is presented in Section~3.  In Section 4 we discuss an extension of Theorem \ref{thm: main} to the case of a positive cosmological constant.  Section 5 includes some further remarks concerning 
Conjecture~\ref{conj.}.

\section{Spacetime and mean curvature flow preliminaries}

We recall some basic definitions and facts. By a spacetime, we mean an $(n+1)$-dimensional, $n \geq 2$, time-orientable Lorentzian manifold $(M,g)$, which, for the analytic purposes of this paper, we assume smoothness of both the manifold and metric. Classically, a spacetime $(M,g)$ is globally hyperbolic if it is strongly causal and all `causal diamonds' $J^+(p) \cap J^-(q)$ are compact. A Cauchy (hyper)surface for a spacetime $(M,g)$ is a set $S$ that is met by every inextendible causal curve exactly once, which implies that $S$ is a topological (in fact locally Lipschitz) hypersurface and $M$ is homeomorphic to $\R \times S$.
It is a well-known fact that a spacetime $(M,g)$ is globally hyperbolic if and only if it contains a Cauchy surface $S$.

Fix a globally hyperbolic spacetime $(M,g)$. Let $F_0 \colon V \to M$ be an embedding such that $V_0 := F_0(V)$ is a smooth compact spacelike Cauchy surface.  The following evolution problem defines mean curvature flow with prescribed constant mean curvature $c \in \R$:
\begin{equation}\label{eq: evol prob}
\frac{d}{ds} F(x,s) \,=\, [(H -c)\nu](x,s) \quad \text{ and } \quad  F(x,0) \,=\, F_0(x),
\end{equation}
where $\nu(x,s)$ is the future-directed unit normal to the smooth spacelike Cauchy surfaces
\begin{equation}\label{eq: V_s}
V_s \,=\, \text{Image}\big[F(\cdot,s)\big],
\end{equation}
and $H(x,s)$ is the mean curvature of $V_s$ at $F(x,s)$ defined with respect to $\nu(x,s)$.  By \cite[Thm. 4.1]{Ecker_Huisken}, a solution to this evolution problem exists at least locally for $s \in [0,s_0)$; that is, there is a smooth map
\begin{equation}
F \colon V \times [0,s_0) \to  M
\end{equation}
such that 
\begin{equation}\label{eq: evol prob2}
F_*(\pd_s) \,=\, (H-c)\nu \quad \text{ and } \quad F(\cdot, 0) \,=\, F_0.
\end{equation}
Moreover, if the flow $\{V_s\}_{s\in [0,s_0)}$ remains in a compact region of the spacetime $M$, then the solution can be extended beyond $s_0$.\footnote{More generally,  
\cite[Thm. 4.1]{Ecker_Huisken} allows not just for a constant $c$ in \eqref{eq: evol prob2}, but for any  smooth function 
$\mathscr{H}$ on $M$ satisfying a certain monotonicity condition.} 

Since $(M,g)$ is globally hyperbolic, we know  from \cite{Sanchez0} that it admits a smooth \emph{temporal function} 
\begin{equation}\label{eq: temporal fun}
t \colon M \to \R.
\end{equation}
That is, $t$ is a smooth function such that its gradient $\nabla t$ is a past-directed timelike vector field and the slices of constant $t$ are smooth spacelike Cauchy surfaces. We fix such a temporal function $t$ on $(M,g)$ and set

\begin{equation}
T \,=\, -\psi \nabla t \quad \text{ where } \quad \psi^{-2} \,=\, - \<\nabla t, \nabla t\>,
\end{equation}
so that $T$ is a future directed  unit timelike vector field.  
For the flow $V_s$, we define the height and gradient
 functions $u,v\colon V \times [0,s_0) \to \R$ as 
\begin{equation}\label{eq: u and v}
u(x,s) \,=\, t\big(F(x,s)\big) \quad \text{ and } \quad v \,=\, -\<\nu, T\>.
\end{equation}
Note that the gradient function $v$ measures the angle between $V_s$ and the slices of constant $t$. 

For the proof of Theorem \ref{thm: main}, we will make essential use of the equation, 
\beq
\frac{d}{ds}u \,=\, (H - c)\psi^{-1}v \label{eq: du/ds},
\eeq
which follows from  \eqref{eq: evol prob} and \eqref{eq: u and v} (cf. 
\cite[Eq. 5]{Ecker_Huisken}).

 

%

\section{Proof of Theorem \ref{thm: main}}

Let $(M,g)$ be a globally hyperbolic spacetime, and consider a (not necessarily smooth) Cauchy surface $S$. Following \cite{AGH}, we say that $S$ has mean curvature $H \leq c$ \emph{in the support sense}, for some number $c$, provided for each $p \in S$ and $\e > 0$
 there exists a smooth spacelike hypersurface $W_{p,\e}$ such that (i) $W_{p,\e}$ is a future support hypersurface for $S$ at $p$, i.e., $W_{p,\e}$ passes through $p$ and is locally to the future\footnote{That is, there is a spacetime neighborhood $U$ of $p$, in which $W_{p,\e}$ is acausal and edgeless, such that $W_{p,\e} \cap U \subset J^+(S,U)$.} of $S$ near $p$, and (ii) the mean curvature of $W_{p,\e}$ at $p$ satisfies $H_{W_{p,\e}}(p) \leq c + \e$. 
A similar definition holds for $S$ having mean curvature $H \geq c$  in the support sense. 

Recall that for each $q \in J^+(S)$, the Lorentzian distance function $d(S,q)$ is defined as the supremum of $L(\g)$ ($=$ length of $\g$) over all causal curves $\g$ from $S$ to $q$.

\begin{prop}\label{prop: small mean curvature}
Let $(M,g)$ be an $(n+1)$-dimensional spacetime with compact Cauchy surfaces. Let $S$ be a fixed Cauchy surface. Consider the sets 
\[
S_\tau = \{q \in M \mid d(S,q) = \tau\}.
\]
Suppose $(M,g)$ is future timelike geodesically complete and satisfies the strong energy condition. Then, for each $\tau > 0$, $S_\tau$ is a Cauchy surface with mean curvature $H \leq \tfrac{n}{\tau}$ in the support sense.
\end{prop}

\proof 
For the convenience of the reader, we provide a complete proof that $S_\tau$ is a Cauchy surface (cf. \cite[Lem. 3.15]{horo1}). We will use the following well-known fact: Since $S$ is compact and $(M,g)$ is globally hyperbolic, the function $x \mapsto d(S,x)$ is continuous, finite-valued, and for each $q \in S_\tau$, there is a unit speed timelike geodesic $\g \colon [0, \tau] \to M$ that realizes the distance from some $p \in S$ to $q$, i.e., $L(\g) = d(S,q)$.

 $S_\tau$ is closed by continuity of the distance function. Moreover, $S_\tau$ is acausal: Clearly $S_\tau$ is achronal. Let $p,q, \g$ be as in the first paragraph of this proof. If there is another point $q' \in S_\tau$ with $q' \in J^+(q)$, then there is a null geodesic  $\l$ from $q$ to $q'$. Consider a convex neighborhood $U$ about $q$. Let $x,y \in U$ be points such that $x \in \g$ and $y \in \l$. The geodesic $\s \subset U$ from $x$ to $y$ has length strictly greater than the portion of $\g$ from $x$ to $q$ \cite[Prop. 5.34]{ON}. Therefore the causal curve formed by following $p$ to $x$ along $\g$, then $x$ to $y$ along $\s$, and lastly $y$ to $q'$ along $\l$ has length strictly greater than $\tau$ -- a contradiction.
Thus $S_\tau$ is closed and acausal. Moreover, it can be seen that $S_\tau$ has no edge points and consequently is a topological hypersurface \cite[Cor. 14.26]{ON}.

Now we show that $S_\tau$ is compact. Let $h$ be a complete Riemannian metric on $M$. If $S_\tau$ is not compact, then there is a sequence of points $q_n \in S_\tau$ such that the $h$-distance from $S$ to $q_n$ tends to infinity. Let $\g_n \colon [0,b_n] \to M$ be an $S$-maximizer from $p_n \in S$ to $q_n$, parameterized by $h$-arclength so that $b_n \to \infty$. Since $S$ is acausal, each $\g_n$ is a timelike geodesic (when appropriately reparameterized). Let $\wt{\g}_n \colon [0, \infty) \to M$ be the future inextendible timelike geodesic extension of $\g_n$, also parameterized by $h$-arclength. By compactness of $S$, there is a $p \in S$ such that $p_n \to p$ (after passing to a subsequence). The limit curve lemma guarantees the existence of a future inextendible causal curve $\g \colon [0, \infty) \to M$ such that $\wt{\g}_n$ limits to $\g$ uniformly on compact subsets of $[0, \infty)$. Moreover, $\g$ is an $S$-maximizer since each $\g_n$ is an $S$-maximizer. (This follows by upper semicontinuity of the Lorentzian length functional.) Therefore, after a reparameterization, $\g$ is a causal geodesic and, in fact, timelike since $S$ is acausal.  By future timelike completeness, there is a $c > 0$ such that $\g(c) \in S_\tau$. Set 
\[
q = \g(c), \quad q' = \g(c+1), \quad q_n = \wt{\g}_n(c), \quad q_n' = \wt{\g}_n(c+1).
\]
Since $b_n \to \infty$, we can choose $n$ large enough so that $c + 1 < b_n$. For such large $n$, we have $q_n, q_n' \in I^-(S_\tau)$, but $q_n' \to q' \in I^+(S_\tau)$ -- yielding an acausality  violation of $S_\tau$. Hence, $S_{\tau}$ is compact.

If $S_\tau$ is not a Cauchy surface, then its Cauchy horizon  is nonempty. Suppose, say, $x \in H^+(S_\tau)$. Then, since $S_\tau$ is closed and acausal, we can apply  \cite[Prop. 14.53]{ON} to conclude that $x$ is the future endpoint of a past-inextendible null geodesic $\eta \subset H^+(S_\tau)$, but this leads to a strong causality violation, 
\[\eta \subset J^-(x) \cap H^+(S_\tau) \subset J^-(x) \cap 	J^+(S_\tau),\] 
since the latter set is compact.


Lastly, we show that $S_\tau$ has mean curvature $H \leq \tfrac{n}{\tau}$ in the support sense. Let $p,q,\gamma$ be as in the first paragraph of this proof.  In particular, each segment of $\g$ is maximal and so there are no cut points to $\g(t)$ along $\g|_{[t, \tau]}$, for $t > 0$. (Recall that a \emph{cut point} to a timelike geodesic segment corresponds to the first point in which the geodesic is no longer maximizing from its starting point. See chapter 9 of \cite{BEE}. This ensures the smoothness of the exponential map near the geodesic, before the cut point, from its starting point. See \cite[Thm. 9.12]{BEE} and \cite[Prop. 10.10]{ON}.)  Consequently, each future Lorentzian sphere $\S_r = \{x \in M \mid d\big(\g(t), x\big) = r\}$, for $r \in (0, \tau - t]$, is smooth near $\g(r+t)$. It follows that $\Sigma_{\tau - t}$ lies locally to the future of $S_\tau$ near $q$, for otherwise there would be a causal curve from $S$ to $S_\tau$ with length greater than $\tau$. Hence $\S_{\tau - t}$ is a future support hypersurface for $S_\tau$ at $q$.
Let $H(r)$ be the mean curvature of $\Sigma_r$ at $\gamma(r+t)$. The Raychaudhuri equation \cite{HE} (for a hypersurface orthogonal geodesic vector field) gives, when restricted to $\g$, 

\beq \label{eq: ray 1}
\frac{dH}{dr} = -\ric(\g', \g') - \frac{H^2}{n} - 2\sigma^2 \leq -\frac{H^2}{n}. 
\eeq
Standard ODE comparison results yield $H(r) \leq \frac{n}{r}$ (see e.g.,  \cite[Section 1.6]{Karcher}), \cite[Section~2]{Esch87}). Evaluating at $r = \tau - t$ and letting $t \to 0$ gives the result.  
\qed

\medskip

\noindent\emph{Remark.}
In fact for $\tau$ small, and provided that $S$ is smooth and spacelike, $S_\tau$ will be a Cauchy surface even without the completeness assumption. More concretely, define the ``future existence time"
 $T_0 > 0$ from $S$ to be the largest value (possibly $ = \infty$)  such that $\exp_p(tN)$ is defined for all $p \in S$ and all $t < T_0$.  (Here $N$ is the future timelike unit normal field to $S$.)  Then, by essentially the same argument, $S_{\t}$ will be a Cauchy surface for all $\t < T_0$.  In this case, for such an $S$,  the conclusion of Proposition \ref{prop: small mean curvature} holds for all $\tau < T_0$.

\medskip

\begin{prop}\label{prop: support barriers}
Let $(M,g)$ be a spacetime with compact Cauchy surfaces satisfying the strong energy condition. Suppose that $S_1$ is a Cauchy surface with mean curvature $H_1 \geq a$ in the support sense and that $S_2 \subset I^+(S_1)$ is a  Cauchy surface with mean curvature $H_2 \leq b$ in the support sense, for some constants $b < c < a$. Then there exists a smooth spacelike Cauchy surface $S$ between\footnote{I.e., $S \subset I^+(S_1) \cap I^-(S_2)$} $S_1$ and $S_2$ with constant mean curvature $H =c$, which is obtained as the asymptotic limit of mean curvature flow.
\end{prop}

\proof
Let $V_0 = F_0(V)$ be any embedded smooth spacelike Cauchy surface between $S_1$ and $S_2$. Consider the evolution problem \eqref{eq: evol prob}. We show that the corresponding flow $V_s$, given by \eqref{eq: V_s}, always lies within the compact region between $S_1$ and $S_2$. Indeed, suppose there is a point $(x_0, s_0) \in V \times [0, \infty)$ where the flow $V_s$ touches $S_2$ for the first time. At this point $p_0 :=F(x_0, s_0)$, we must have  
$\tfrac{d}{ds}u \geq 0$, as otherwise the flow $V_s$ would necessarily 
have touched $S_2$ at some value $s_1 < s_0$.
 We conclude from \eqref{eq: du/ds} that $H(x_0, s_0) \geq c$. Since $S_2$ has mean curvature $H_2 \leq b < c$ in the support sense, for each $\e > 0$, there exists a future support hypersurface $W$ for $S_2$ at $p_0$ with mean curvature $H_W \leq b + \e$. By choosing $\e$ sufficiently small, it is easy to see that there is a $d \in (b,c)$ such that $H_{V_{s_0}} > d$ at $p_0$ and $H_{W} < d$ at $p_0$. Since $W$ is locally to the future of $S_2$ near $p_0$ and hence also locally to the future of $V_{s_0}$ near $p_0$, we arrive at a contradiction to the 
  maximum principle. Thus the flow always remains below $S_2$. Likewise, the flow always remains above $S_1$. 

Since the Cauchy surfaces are compact, we can find $t$-slices of the temporal function \ref{eq: temporal fun} such that $S_1$ always lies above a slice $t = \a$ and $S_2$ lies below a slice $t= \b$. Therefore the flow always remains in the smooth compact region of spacetime between $t = \a$ and $t= \b$. It follows via \cite[Thms. 4.1 and 4.2]{Ecker_Huisken} that the flow $V_s$ is defined for all $s \in [0, \infty)$ and (by taking a suitable subsquence) limits to a smooth spacelike Cauchy surface $S$ with constant mean curvature $H = c$.
\qed

\medskip

\noindent\emph{Remark.} As alluded to in the introduction, Proposition \ref{prop: support barriers} extends Theorem 4.2 in \cite{Ecker_Huisken} to the case of mean curvature barriers in the support sense.
An existence result similar to \cite[Thm. 4.2]{Ecker_Huisken}
was  obtained by Claus Gerhardt \cite{Gerhardt}, again assuming smooth barriers. Gerhardt's proof was via a variational method and did not require the Ricci curvature condition.   There appears to be certain technical issues in trying to extend Gerhardt's proof to the case of barriers in the support sense.  A closely related situation arises in \cite{ABBZ}. We thank the authors of \cite{ABBZ} for a discussion pertaining to some of these issues.

%

\medskip

\noindent\emph{Proof of Theorem \ref{thm: main}}. By assumption, there is a smooth spacelike Cauchy surface $V$ with mean curvature $H_V \geq 0$. If $H_V = 0$ identically, then we are done. Otherwise, there is a point $p \in V$ with $H_V(p) > 0$. We can perturb $V$ slightly to the past to produce a smooth spacelike Cauchy surface with strictly positive mean curvature. To see this, let $V_r$ be a variation of $V_0 = V$ with variation vector $\phi \nu$, where $\phi$ is a smooth function on $V$ and $\nu$ is the future directed unit timelike normal to $V$. 
 Let $H_r$ be the mean curvature of $V_r$. A standard computation shows that 
\[
\frac{\pd H_r}{\pd r}\bigg|_{r=0} \,=\, -L\phi \quad \text{ where } \quad L\phi \,:=\, -\Delta \phi + \big(\ric(\nu,\nu) + \tfrac{H_0^2}{n} + \sigma^2 \big)\phi,
\]
where $\Delta$ is the spatial Laplacian on $V$ and $\sigma^2$ is the norm squared of the trace-free part of the second fundamental form of $V$. Since $H_0(p) > 0$, the Rayleigh quotient formula for the principal eigenvalue $\l_1$ of $L$,
\[
\lambda_1 \,=\, \inf_{1 = ||f||_{L^2(V)} }  \int_V\left(|\nabla f|^2 + \big(\ric(\nu,\nu) + \tfrac{H_0^2}{n} + \sigma^2\big)f^2\right),
\]
 shows that $\lambda_1 > 0$. Then, by varying with respect to the principal eigenfunction, there are times $r < 0$ such that $H_r > 0$ everywhere on $V_r$.

 Thus we can assume -- without loss of generality -- that $H_V > 0$. Since $V$ is smooth, we can choose constants $a > c > 0$ such that $H_V \geq a$. 
 By taking $\tau$ sufficiently large, Proposition \ref{prop: small mean curvature}  implies the existence of a Cauchy surface $V_\tau$ with mean curvature $H_\tau \leq b$ in the support sense, for some $b < c$. Then Proposition \ref{prop: support barriers} implies that there is a smooth spacelike Cauchy surface between $V$ and $V_\tau$ with constant mean curvature $H = c$. \qed
 
\medskip


\medskip
\noindent
\emph{Remark.} The timelike future completeness assumption in Theorem \ref{thm: main} can be weakened some. Suppose, for instance, $V$ has mean curvature $H_V > c$, then, instead of future completeness, it is sufficient that the future existence time $T_0$ introduced in the remark after the proof of
Proposition \ref{prop: small mean curvature} satisfy $T_0 >  \frac{n}{c}$. 
In particular, in a big bang scenario, we imagine $c$ being large (e.g., a crushing singularity) in which case we would only require a small $T_0$ to ensure a CMC slice.

\medskip

\noindent\emph{Remark.} In Theorem \ref{thm: main}, we showed the existence of a smooth Cauchy surface, say $V_1$, with constant mean curvature $H = c_1 > 0$ (in the case that $V$ is not maximal); moreover, $V_1$ is unique since $c_1 \neq 0$, see e.g., \cite[Thm. 7.3]{Gerhardt}. In fact, we can say more. By another application of Theorem \ref{thm: main}, for any $0 < c_2 < c_1$, we can find a CMC Cauchy surface $V_2 \subset I^+(V_1)$ with mean curvature $H = c_2$. Then the region between $V_1$ and $V_2$ is covered by unique CMC Cauchy surfaces with mean curvatures $c \in (c_2, c_1)$, as follows from \cite[Thm. 7.5]{Gerhardt}.

\section{The case of a positive cosmological constant}

In this section we indicate how Theorem \ref{thm: main} extends to the case of a positive cosmological constant.  

\subsection{Extending Theorem \ref{thm: main} to the case $\Lambda > 0$}

The Einstein equation with a cosmological constant $\Lambda$ when solved for the spacetime Ricci tensor becomes
\beq\label{eq: ric}
\text{Ric} \,=\, 8\pi\left(T - \frac{1}{n-1}\text{tr}_g(T) g\right) + \frac{2}{n-1}\Lambda g.
\eeq
By imposing a standard condition on the energy-momentum tensor \cite[Section 4.3]{HE},
one is led to the following version of the strong energy condition in the case of a positive cosmological constant,
\beq\label{eq: lambda-SEC}
\ric(X,X) \ge -n\l \,  \quad \text{for all unit timelike vectors  }  X  \,,
\eeq
where $\l = \frac{2}{n(n-1)}\Lambda > 0$.

In the case of a positive cosmological constant, we have the following analogue of Proposition \ref{prop: small mean curvature}.

\begin{prop}\label{prop: small mean curvature 2}
Let $(M,g)$ be an $(n+1)$-dimensional spacetime with compact Cauchy surfaces. Fix a Cauchy surface $S$, and consider the sets 
\[
S_\tau = \{q \in M \mid d(S,q) = \tau\}.
\]
Suppose $(M,g)$ is future timelike geodesically complete and satisfies \eqref{eq: lambda-SEC}.  Then, for each $\tau > 0$, $S_\tau$ is a Cauchy surface with mean curvature $H \leq 
n\sqrt{\l}\coth(\sqrt{\l}\tau)$ in the support sense.
\end{prop}

\proof The proof proceeds just as in the proof of Proposition \ref{prop: small mean curvature}.
In  the present situation the inequality \eqref{eq: ray 1} becomes,
$$
\frac{d\mathcal{H}}{dr} \leq \lambda - \mathcal{H}^2,
$$
where $\mathcal{H} =  \frac{H}{n}$.   The theorem again follows from standard comparison results (see e.g.,  \cite[Section 1.6]{Karcher}), \cite[Section~2]{Esch87}). \qed

\medskip 

As discussed in the next subsection, Theorem 4.1 and the existence part of Theorem 4.2 in \cite{Ecker_Huisken} remain valid if the strong energy condition is replaced by the inequality \eqref{eq: lambda-SEC}. 
We thus obtain the following extension of Theorem 
\ref{thm: main}.

\begin{thm}\label{thm: main2}
Let $(M,g)$ be a spacetime with compact Cauchy surfaces, and suppose there exists a smooth spacelike Cauchy surface with mean curvature $H \ge  n\sqrt{\l}$.  If $(M,g)$ is future timelike geodesically complete and the Ricci tensor satisfies the inequality \eqref{eq: lambda-SEC}, then $(M,g)$ contains a CMC Cauchy surface.
\end{thm}

\subsection{Remarks on extending the results in \cite{Ecker_Huisken} to the case  $\Lambda > 0$}\label{appendix}

One of the assumptions in \cite{Ecker_Huisken} is that the spacetime $(M,g)$ satisfies the strong energy condition, i.e., $\ric(X,X) \geq 0$ for all timelike $X$. As discussed above, when the Einstein equations hold with a positive cosmological constant, the appropriate version of this energy condition is the assumption \eqref{eq: lambda-SEC}, where $\lambda > 0$.

The work in \cite{Ecker93} in the cosmological setting shows that the existence results in \cite{Ecker_Huisken} can be extended to the case when the spacetime satisfies \eqref{eq: lambda-SEC}. Moreover, the forcing term $\mathscr{H}$ in \cite{Ecker93} can be any smooth function; e.g., it need not be constant as it is in this work nor, more generally, satisfy the monotonicity assumption as in \cite{Ecker_Huisken}.  In this subsection, we wish to point out  how the arguments in \cite{Ecker_Huisken} can be easily modified to handle the case of a positive cosmological constant in our setting.

The proof of Theorem 4.1 in \cite{Ecker_Huisken} (and hence also the proof of the existence part of Theorem 4.2 in \cite{Ecker_Huisken}) relies on the estimates given in Propositions 4.4, 4.6, and 4.7 in \cite{Ecker_Huisken}. Specifically, Propositions 4.4 and 4.7 rely on the mean curvature estimate established in 4.6 -- this is the only place where the energy condition plays a role. 

We show how \cite[Prop. 4.6]{Ecker_Huisken} generalizes to the case when \eqref{eq: lambda-SEC} is assumed. However, in this new setting, the constant $C$ appearing in the statement of \cite[Prop. 4.6]{Ecker_Huisken} will now also depend on $\lambda$. From Proposition 3.3(ii) and eq. (1) in \cite{Ecker_Huisken}, we have the following equation 
\beq\label{eq: flow of H2}
\left(\frac{d}{ds} - \Delta_s \right)(H - \mathscr{H}) \,=\, - (H- \mathscr{H})\big(|A|^2 + \Ric(\nu, \nu) + \< \nabla\mathscr{H}, \nu \> \big),
\eeq
where $\mathscr{H}$ is a smooth function on $M$ and $\Delta_s$ is the spatial Laplacian induced from the Riemannian metric on $V_s$. (Note there is a typo in the proof of \cite[Prop. 4.6]{Ecker_Huisken}. A parenthesis is misplaced in the first line of the proof.)

Set $f = H - \mathscr{H}$. Then
\[
\left(\frac{d}{ds} - \Delta_s \right)f^2
 \,=\, -2f^2 \big(|A|^2 + \Ric(\nu, \nu) + \< \nabla\mathscr{H}, \nu \> \big) - 2|\nabla_s f|^2.
\]
Assuming \eqref{eq: lambda-SEC} and  the monotonicity assumption as in \cite{Ecker_Huisken}, $\< \nabla\mathscr{H}, \nu \> \geq 0$ (which is trivially satisfied in our situation $\mathscr{H}=c=$ const.), we have
\begin{align*}
\left(\frac{d}{ds} - \Delta_s \right)f^2 \,&\leq\, -2f^2|A|^2 + 2n\lambda f^2 
\\
&\leq\, -\frac{2}{n}f^2H^2 +  2n\lambda f^2 \,.
\end{align*}
Also,
\[
-f^2H^2 \,=\, -f^2(H- \mathscr{H} + \mathscr{H})^2 \,\leq\, -f^4 + 2|f^3\mathscr{H}| \,.
\]
Using the Peter-Paul inequality ($ab \leq \tfrac{a^2}{2\e} + \tfrac{\e  b^2}{2}$) twice, we estimate the cubic term:
\begin{align*}
|f^3\mathscr{H}| \,&\leq\, f^2\left(\frac{f^2}{2\e_1} + \frac{\e_1\mathscr{H}^2}{2} \right)
\\
&\leq\, \frac{f^4}{2\e_1} + \frac{\e_1}{2}\left(\frac{f^4}{2\e_2} + \frac{\e_2\mathscr{H}^4}{2} \right).
\end{align*}
Therefore
\[
\left(\frac{d}{ds} - \Delta_s \right) f^2 \,\leq\, -\frac{2}{n}f^4 + \frac{4}{n}\left(\frac{1}{2\e_1} + \frac{\e_1}{4\e_2} \right)f^4 + n \e_1\e_2\mathscr{H}^4 + 2n\lambda f^2.
\]
Using Peter-Paul again, we have
\[
2n\lambda f^2 \,\leq\, 2n\lambda\left(  \frac{f^4}{2\e_3} + \frac{\e_3}{2}\right).
\]
Therefore
\[
\left(\frac{d}{ds} - \Delta_s \right) f^2 \,\leq\, -\frac{2}{n}f^4 + \frac{4}{n}\left(\frac{1}{2\e_1} + \frac{\e_1}{4\e_2} \right)f^4 + n \e_1\e_2\mathscr{H}^4 + \frac{n\lambda}{\e_3}f^4 + n\lambda \e_3.
\]
We can choose $\e_1, \e_2, \e_3 > 0$ such that the coefficient in front of $f^4$ is $\leq -\tfrac{1}{n}$. Hence
\[
\left(\frac{d}{ds} - \Delta_s \right) f^2 \,\leq\, -\frac{1}{n}f^4+ n \e_1\e_2\mathscr{H}^4 + n\lambda \e_3.
\]
This inequality is equivalent to the last inequality in the proof of \cite[Prop. 4.6]{Ecker_Huisken}, where the constant now also depends on $\Lambda$.
The rest of the proof is as in \cite{Ecker_Huisken}.

\section{ Further comments on Conjecture \ref{conj.}}


In this last section, we once again assume the strong energy condition. We provide some further comments on Conjecture \ref{conj.} and its relationship to the future causal boundary. 
In section \ref{sec: spacelike causal bdy} we discuss some issues concerning the structure of the future causal boundary of cosmological spacetimes.  In 
section~\ref{sec: maximal cauchy} we establish a connection between the future causal boundary being spacelike and the existence of a maximal Cauchy surface. The relationship to the Bartnik splitting conjecture is also discussed.

\subsection{Some remarks on the causal boundary of cosmological spacetimes}\label{sec: spacelike causal bdy}

Conjecture \ref{conj.} was also motivated in part by the following conjecture \cite{Bart88}.   

\medskip

\begin{conj}[Bartnik splitting conjecture]\label{conj2.}
Let $(M,g)$ be a spacetime which contains a
compact Cauchy surface and obeys the strong energy condition,
$\ric(X,X)\ge0$ for all timelike vectors $X$. If $(M,g)$ is timelike geodesically complete,
then $(M,g)$ splits isometrically into the product
$(\bbR \times V,-dt^2\oplus h)$, where $(V,h)$ is a compact Riemannian manifold.
 \end{conj}

For background and motivation,  see e.g., \cite{GalBart}.  
While Conjecture \ref{conj2.} has been proven subject to various auxilliary conditions,
it remains open in full generality.  
However, it is well known that  Conjecture \ref{conj2.} holds if one can establish the existence of a CMC Cauchy surface.  
In particular, it follows from
Theorem \ref{CMCnonpos}  that Conjecture \ref{conj2.} holds under the stronger assumption of nonpositive timelike sectional curvatures.  

We make a few remarks about the proof of Theorem \ref{CMCnonpos}.  The proof involves the notion of the future causal boundary of a spacetime.  Let $(M,g)$ be a globally hyperbolic spacetime. Very briefly, the future causal boundary $\mathscr{C}^+$ consists of the set of all {\it terminal indecomposable past} sets, or TIPs for short. These are the sets of the form $I^-(\g)$, where $\g$ is a future inextendible timelike curve in $M$, with appropriate identifications, when two such curves have the same timelike past.  
See e.g., \cite[Section~6.8]{HE}, for a more detailed discussion. 

Note that $\mathscr{C}^+$ consists of a single element if and only if $I^-(\g) = M$ for all future inextendible timelike curves, i.e., if and only if there are no future observer horizons.  The timelike sectional curvature condition in Theorem \ref{CMCnonpos} is used to show that $\mathscr{C}^+$ consists of a single element.  This then guarantees the existence of a CMC Cauchy surface:  from an observation of Tipler \cite{Tipler}, it follows that a criterion of  Bartnik \cite{Bart88} for the existence of a CMC Cauchy surface is satisfied.  

The authors have previously raised the question \cite{GalBart,GalLing2} as to whether one can show that $\mathscr{C}^+$ consists of a single element if one replaces the timelike sectional curvature assumption  in Theorem \ref{CMCnonpos} with the strong energy condition.  This of course would settle Conjecture \ref{conj2.}.  In  fact, the answer is ``no", as we now illustrate with a specific example.     Consider the multiply warped product spacetime of the form

\beq\label{eq: multi}
M = (0, \infty) \times T^3, \quad g = -dt^2 + \sum_{i = 1}^3 a_i^2(t) (dx^i)^2 \,,
\eeq
where the $a_i$'s are positive, and the $x^i$'s should be viewed as periodic coordinates. (Such spacetimes are sometimes used as spatially homogeneous models in cosmology~\cite{HE}.)  It is convenient to assume each $x^i$ lies in an interval, $- b \le x^i \le b$, with the end points identified.
Basic causal theoretic arguments show that the tori, $V_t = \{t\} \times T^3$ are Cauchy surfaces, and hence that $(M,g)$ is globally hyperbolic.

Now consider the following choice for the $a_i$'s:
\beq\label{eq: choice}
a_1 = a_2 = t^{\frac34} \,, \;   \; a_3 = t^{\frac54} \, .
\eeq
It is not difficult to show by direct computation that the resulting model satisfies the strong energy condition, in fact it satisfies $\ric(X,X) \ge 0$ for all vectors $X$.  (Alternatively, one could apply \cite[Prop. 3.2]{Unal05}.)  As a hint to this, observe that,
\beq\label{r00}
\ric(e_0,e_0) =  - \left(\frac{a_1''}{a_1} + \frac{a_2''}{a_2} + \frac{a_3''}{a_3}\right)
=  - \left(-\frac{3}{16}\frac{1}{t^2} -\frac{3}{16}\frac{1}{t^2}+ \frac{5}{16}\frac{1}{t^2}\right) > 0 \,,
\eeq
where $e_0 = \frac\d{\d t}$.  Moreover,  it is not difficult to show that $M$ is future timelike geodesically complete; this essentially relies on the fact that 
$\int_{t_0}^{\infty} a_i dt = \infty$, $i = 1,2, 3$.  (For a general completeness result applicable to  multiply warped product spacetimes, see \cite[Section 4]{Unal00}.)

Now we consider the structure of the future causal boundary $\mathscr{C}^+$.
Let $t \to \g(t)$ be the $t$-line centered at $x^i = 0$, $i = 1,2,3$.  Let $s \to \eta(s)$ be a past directed null geodesic starting at $(t_0,0,0,0)$ on $\g$, with $t_0 > 1$. With the choice $a_3 = t^{\frac54}$,
it is not difficult to show that the intersection of every such null geodesic $\eta$  with the toroidal  time slice $V_1 = \{t =1\}$ is confined to the region  $-4 < x^3 < 4$.  This boils down to the fact that $\int_1^{t_0} t^{-\frac54} dt  < 4$.   Since $\d I^-(\g(t_0))$ is ruled by such geodesics, it follows that  $I^-(\g(t_0))$ does not meet the region $|x^3| \ge 4$.  Since $I^-(\g) = \cup_{t_0 \in (1, \infty)} I^-(\g(t_0))$, we conclude that  $I^-(\g) \cap V_1$ does not meet the region $|x^3| \ge 4$.  In particular $I^-(\g) \ne M$.  
This immediately implies that $\mathscr{C}^+$ consists of more than one element. However, one can say somewhat more.

From the $x^3$-translational symmetry,  a similar conclusion to the above holds for each $t$-line centered at $x^1 = x^2 = 0, x^3 = \xi$, $-b \le \xi \le b$ ($b >4$), and, moreover, each such $t$-line  defines a distinct element of $\mathscr{C}^+$.  It can be further shown that 
$\mathscr{C}^+$ consists only of these elements.  In fact, much of this follows from a general result of Harris \cite[Prop.~3.5]{Harrisdiscgp}, the relevant part of which we now paraphase. 

By a general multiply warped product spacetime  we mean a spacetime $(M,g)$ of the form,
\begin{align}
M &= (a,b) \times V_1 \times \cdots \times V_m  \,,   \text{ and} \notag \\
g &= -dt^2 + a_1^2(t) h_1 + a_2^2(t) h_2 + \cdots + a_m^2(t) h_m  \,, \notag
\end{align}
where $(a,b)$ is an interval (possibly $a = -\infty$ or $b = \infty$), and where
 for $i = 1,..., m$,
$(V_i,h_i)$ is a  Riemannian manifold.

The following is immediate from \cite[Prop. 3.5]{Harrisdiscgp}.
\begin{prop}
Let $(M,g)$ be a multiply warped product spacetime as above, such that each $V_i$ is compact.  Suppose
\ben
\item[(i)] for $i = 1, ..., k$, and some $c \in (a,b)$, $\int_c^b a_i^{-1} dt = \infty$, and  
\item[(ii)] for $i = k+1, ..., k_m$, and some $c \in (a,b)$, $\int_c^b a_i^{-1} dt < \infty$.
\een
Then the future causal boundary $\mathscr{C}^+$ (in an appropriate topology) is homeomorphic to $V_{k+1} \times \cdots \times V_m$.  Moreover,  $\mathscr{C}^+$ is spacelike.  
\end{prop}

\smallskip
This result applied to our example \eqref{eq: multi}, \eqref{eq: choice} implies that 
$\mathscr{C}^+$ is homeomorphic to $S^1$, as is somewhat suggested by the discussion above.  Furthermore, $\mathscr{C}^+$ is spacelike.
In general, the future causal boundary being \emph{spacelike} may be taken to mean that no TIP is properly contained in another TIP.

Thus, the simple example presented here, which can be generalized in various directions, shows that, under the assumptions of Conjecture \ref{conj.}, the future causal boundary need not consist of a single element.   One, however, may be inclined to formulate a new question:  Under the assumptions of Conjecture \ref{conj.}, will the future causal boundary always be spacelike? An example showing otherwise could not, of course, be within this class of multiply warped product spacetimes.  Moreover,  an affirmative answer would settle the Bartnik splitting conjecture. For, as shown in \cite{horo2}, in this case the spacetime admits a timelike line (i.e., a globally maximizing inextendible timelike geodesic), and one can apply the Lorentzian splitting theorem under the assumptions of the Bartnik splitting conjecture.  In the next section we provide an alternative proof of this fact. We show that if a spacetime satisfies the hypotheses of Conjecture \ref{conj.} but is also past timelike complete 
and has a spacelike future causal boundary, then the spacetime contains a maximal  Cauchy surface and consequently splits \cite{Bart88}.

\subsection{Existence of a maximal Cauchy surface from a spacelike causal boundary}\label{sec: maximal cauchy}

Let $(M,g)$ be a globally hyperbolic spacetime and let $S$ be a (not necessarily smooth) Cauchy surface. A future inextendible timelike geodesic $\g \colon [0, a) \to M$, parameterized by proper time, is called a future $S$-ray provided it maximizes the distance to $S$, i.e., $L(\g|_{[0,t]}) = d\big(S, \g(t)\big)$ for all $t \in [0,a)$. A past $S$-ray is defined similarly. If $S$ is compact, then it always admits a future and past $S$-ray. If, for some Cauchy surface $S$,  $M$ contains a future $S$-ray $\g$ and a past $S$-ray $\eta$ such that $I^-(\g) \cap I^+(\eta) \neq \emptyset$, then we say that ``$M$ satisfies the ray-to-ray condition."

\begin{thm}\label{thm: maximal}
Suppose $(M,g)$ satisfies the assumptions of Conjecture \ref{conj.} but is also past timelike complete. If $M$ satisfies the ray-to-ray condition, then $M$ contains a maximal (i.e. mean curvature zero) Cauchy surface. 
\end{thm}

\smallskip

A spacelike future causal boundary implies the ray-to-ray condition, see  \cite[Prop.~5.9]{horo2}. Thus we have the following corollary.

\begin{cor}
Suppose $(M,g)$ satisfies the assumptions of Conjecture \ref{conj.} but is also past timelike complete. If $M$ contains a spacelike future causal boundary, then $M$ contains a maximal Cauchy surface. 
\end{cor}

\smallskip

\noindent\emph{Proof of Theorem \ref{thm: maximal}}. For some Cauchy surface $S$, there is  a future $S$-ray, $\g \colon [0, \infty) \to M$, and a past $S$-ray, $\eta \colon [0, \infty) \to M$, such that $I^-(\gamma) \cap I^+(\eta) \neq \emptyset$. Therefore there are points $p = \eta(t_1)$ and $q = \gamma(t_2)$ such that $q \in I^+(p)$. Set
\[
S_1 \,=\, \{x \in M \mid d(x,S) = t_1\} \quad \text{ and } \quad S_2 \,=\, \{x \in M \mid d(S,x) = t_2\}.
\]
$S_1$ and $S_2$ are Cauchy surfaces by Proposition \ref{prop: small mean curvature}. Also, as is easily seen, $\g|_{[t_2, \infty)}$ is a future $S_2$-ray and $\eta|_{[t_1, \infty)}$ is a past $S_1$-ray. 

Let $S^-_\infty$ be the \emph{past ray horosphere} \bl{(\cite{horo1})} associated to $\g|_{[t_2, \infty)}$, that is, 
\[
S^-_\infty \,=\, \pd \left(\bigcup_{k = 1}^\infty I^-(S^-_k) \right),
\]
where $S^-_k$ are the Lorentzian spheres
\[
S^-_k = \{x \mid d\big(x, \g(t_2 + k) \big) =k\}.
\]
$S^-_\infty$ enjoys the following nice properties:
\begin{itemize}
\item[(1)] $S^-_\infty$ is nonempty, in particular $q \in S^-_\infty$.
\item[(2)] $S^-_\infty$ is acausal and edgeless.
\item[(3)] $S^-_\infty \subset J^-(S_2)$ and is a past Cauchy surface, i.e., $H^-(S^-_\infty) = \emptyset$. 
\item[(4)] $S^-_\infty$ has mean curvature $H \geq 0$ in the support sense (with one-sided Hessian bounds).
\end{itemize}

We briefly allude to the proof of these properties. $S^-_\infty$ is an achronal boundary by construction and hence it is edgeless. Moreover, $S^-_\infty$ does not intersect $I^+(S_2)$ and so $S^-_\infty \subset J^-(S_2)$; consequently, $S^-_\infty$ is a past Cauchy surface by \cite[Prop. 3.17]{horo1} and acausal by \cite[Thm. 3.18]{horo1}. Lastly, (4) follows from a similar argument as in Proposition \ref{prop: small mean curvature}, see \cite[Thm. 4.2]{horo1}.

Similarly, let $S^+_\infty$ be the \emph{future ray horosphere} associated to $\eta|_{[t_1, \infty)}$. It also satisfies analogous properties of (1) - (4) above. Consider the compact set $K := J^+(S_1) \cap J^-(S_2)$. Using (3) above and the analogous property for $S^+_\infty$, 
it is not hard to see that 
\[
J^-(S^-_{\infty}) \cap J^+(S^+_{\infty}) \,=\, J^-(S^-_\infty \cap K) \cap J^+(S^+_\infty\cap K).
\]
Consequently, $J^-(S^-_{\infty}) \cap J^+(S^+_{\infty})$ is compact. The existence of a maximal Cauchy surface now follows from \cite[Thm. 3.7]{GalBanach}. The proof of this theorem is somewhat involved, so we give a more direct proof.

Since $J^+(S^+_\infty) \cap J^-(S^-_\infty)$ is compact, there exist points $x_1 \in S^+_\infty$ and $x_2 \in S^-_\infty$ such that $\tau := d(x_1, x_2)$ realizes the Lorentzian distance between $S^+_\infty$ and $S^-_\infty$. Now set $\Sigma = \{x \in M \mid d(x,S^+_\infty) = \tau\}$; using the fact that $S^+_\infty$ is a future Cauchy surface, arguments as in the proof of Proposition \ref{prop: small mean curvature} show  that $\Sigma$ is a closed acausal $C^0$ hypersurface, and moreover it lies locally to the future of $S^-_\infty$ near $x_2$. By \cite[Lem. 3.6]{GalBanach} $\Sigma$ has mean curvature $H \leq 0$ in the support sense. (The proof involves, modulo a certain technicality, pushing the support hypersurfaces of $S^+_\infty$ along the normal geodesics and making use of the Raychaudhuri equation.)

Hence to summarize, there is a spacetime neighborhood $U$ of $x_2$ such that $\Sigma$ lies to the future of $S^-_\infty$ within $U$ and $\Sigma$ has mean curvature $H \leq 0$ in the support sense, whereas $S^-_\infty$ has mean curvature $H \geq 0$ in the support sense. It then follows from a rough version of the Lorentzian geometric maximum principle, \cite[Thm. 3.6]{AGH}, that there is an open set $A \subset \Sigma \cap S^-_\infty$ which is a smooth spacelike hypersurface with mean curvature $H = 0$.\footnote{Proposition 3.5 in \cite{AGH} guarantees that, in our situation, a certain technical condition in the statement of \cite[Thm.~3.6]{AGH} is satisfied.} A continuity argument then shows that this intersection globalizes. Thus $\Sigma = S^-_\infty$ is a smooth maximal hypersurface which is closed and acausal. Since it is contained in the compact set between $S_1$ and $S_2$, it is also compact and consequently a Cauchy surface.
\qed

\section*{Acknowledgments} 
We thank Piotr Chru{\'s}ciel, Stacey Harris, Gerhard Huisken, and  Alex Mramor for helpful comments and discussions. Gregory Galloway was supported by  the Simons Foundation, Award No. 850541. Eric Ling was supported by Carlsberg Foundation CF21-0680 and Danmarks Grundforskningsfond CPH-GEOTOP-DNRF151.  Part of the research on this paper was supported by the National Science Foundation under Grant No.
DMS-1928930 while the authors were in residence at the Simons Laufer Mathematical Sciences Institute (formerly MSRI) in Berkeley, California, during the Fall 2024 semester.

\vspace{.1in}
\noindent
{\bf Declarations.} On behalf of both authors, the corresponding author states that there is no conflict of interest. No data was collected or analysed as part of this project.

\newpage

\bibliographystyle{amsplain}

\end{document}